# On a Theorem on sums of the form

$$1+2^{2n}+2^{2n+1}+\ldots\ldots+2^{2n+m}$$

# and a result linking Fermat with

# Mersenne numbers

*By Konstantine 'Hermes' Zelator*

*Department Of Mathematics*

*College Of Arts And Sciences*

*Mail Stop 942*

*University Of Toledo*

*Toledo, OH 43606-3390*

*U.S.A.*

# 1. INTRODUCTION

In this paper we study sums of the form $1 + 2^{2^n} + 2^{2^{n+1}} + \ldots + 2^{2^{n+m}}$.

In his book *250 Problems in Elementary Number Theory*, (see reference [1], problem 123, p. 66), W. Sierpinski shows that the numbers $1 + 2^{2^n} + 2^{2^{n+1}}$ are divisible by 21, for $n = 1, 2, \ldots$. This leads to the fact that for

$$n \geq 2, \frac{1}{3}\left(1 + 2^{2^n} + 2^{2^{n+1}}\right)$$ are all composite numbers divisible by 7 (for $n = 1$, one obtains 7). This problem of W. Sierspinski has motivated the study in this paper of integers of the form $1 + 2^{2^n} + 2^{2^{n+1}} + \ldots + 2^{2^{n+m}}$.

The main theorem in this paper is stated as follows. Let $N$ be an odd integer $\geq 3$, $n$ an integer $\geq 2$ such that $2^n \equiv 1 \bmod s$, where $s$ is the exponent to which 2 belongs modulo $N$ (note that $s$ must be odd). Suppose that $a$ is the exponent to which 2 belongs modulo $s$. Also, let $l$ be the largest integer so that $2^l < s$; thus $2^{l+1} \geq s$ and so $2^{l+1} \equiv b \bmod s$, for some $b$ with $0 \leq b < s$. Consider the integers $k_l = 1 + 2^{2^n} + 2^{2^{n+1}} + \ldots + 2^{2^{n+l}}$ and

$k_b = 2^b + 2^{b.2} + 2^{b.2^2} + \ldots + 2^{b.2^{a-1}}$.

Assume that $d \mid k_l$, where $d = (k_b, N)$. Then there are infinitely many integers $m$ for which

$I_{n,m} = 1 + 2^{2^n} + 2^{2^{n+1}} + \ldots + 2^{2^{n+m}} \equiv 0 \bmod N$.

In fact the integers $m = l + r.a + i.L.D$ have precisely this property; where $L$ is the least common multiple of $N$ and $a$, $D = (N, a)$ and $r$ is the least positive integer that satisfies the congruence

$$k_l + r.k_b \equiv 0 \mod N.$$

(The existence of such, an $r$ is guaranteed by the condition $d \mid k_l$; $r$ can be effectively computed.)

Below Theorem 1, two examples are presented. The first involves $N = 31, n = 4k$ ($k$ any positive integer) and the second $N = (31)(7) = 217, n = 4k$. In this introduction we present an additional example. Take

$$N = 7.127 = 889, n = 6k, \text{ for } k = 1, 2, \ldots.$$

The exponent to which 2 belongs modulo $N$ is $3.7 = 21$ (note that exponent of 2 mod 7 is 3, while the exponent of 2 mod 127 is 7; thus the exponent of 2 mod 7.127 is $3.7 = 21$, the least common multiple of 3 and 7). So $s = 21$, and $2^n = (2^6)^k \equiv 1 \mod 21$. Furthermore, $a = 6; a - 1 = 5$, and l=4; so $2^{l+1} + 2^5 \equiv 11 \mod 21$

which means $b = 11$. We have $k_l = 1 + 2^2 + 2^{2^2} + 2^{2^3} + 2^{2^4}$ and

$$k_b = 2^{11} + 2^{11.2} + 2^{11.2^2} + 2^{11.2^3} + 2^{11.2^4} + 2^{11.2^5}.$$

Let us consider $k_b$ modulo 7 and 127. We have

$$2^{11} \equiv 2^2, 2^{11.2} \equiv 2, 2^{11.2^2} \equiv 2^2$$

$$2^{11.2^3} \equiv 2, 2^{11.2^4} \equiv 2^2, 2^{11.2^5} \equiv 2 \mod 7$$

and so $k_D$ $k_b = 2^2 + 2 + 2^2 + 2 + 2^2 + 2 \equiv 18 \equiv 4 \mod 7$.

On the other hand,

$$2^{11} \equiv 2^4, 2^{11.2} \equiv 2, 2^{11.2^2} \equiv 2^2$$

$$2^{11.2^3} \equiv 2^4, 2^{11.2^4} \equiv 2, 2^{11.2^5} \equiv 2^2 \mod 127$$

Consequently, $k_b = 2^4 + 2 + 2^2 + 2^4 + 2 + 2^2 \equiv 44 \mod 127$.

Thus, since $k_b \equiv 4 \bmod 7$ and $k_b \equiv 44 \bmod 7$, we have $(k_b, 7.127) = 1 = d$ and, trivially $d \mid k_l$, and so the hypothesis of the theorem is satisfied.

Let us now find the least positive integer $r$ so that

$$k_l + r.k_b \equiv \bmod(7.127)$$

Since $k_l \equiv 1 \bmod 7$ and $k_l \equiv 27 \bmod 127$, the last congruence yields

$$1 + 4r \equiv 0 \bmod 7$$

$$27 + 4r \equiv 0 \bmod 127.$$

The first congruence implies $r \equiv 5 \bmod 7$. Write the second as

$27 + 11.2^2.r \equiv 0 \bmod 127$. If we multiply the congruence with $11^{62}.2^5$ (we have $11^{63} \equiv 1, 2^7 \equiv 1 \bmod 127$), we obtain $27(11^{62}.2^5) + r \equiv 0 \bmod 127$ or

$$r \equiv -27.11^{62}.2^5 \equiv 60 \bmod 127.$$

The least positive integer $r$ that satisfies both congruences $r \equiv 5 \bmod 7$ and $r \equiv 60 \bmod 127$, is $r = 1076$.

Also the least common multiple L of $N = 889$ and $a = 6$ is $L = 889.6 = 5334$, while $D = (N, a) = 1$. According to the (main) theorem, it follows, since $m = l + r.a + i.L.D = 4 + 1076.6 + 5334i = 6460 + 5334i$, that the numbers

$$1 + 2^{2^{6k}} + 2^{2^{6k+1}} + \ldots + 2^{2^{6k+(6460+5334i)}} \text{ are divisible by } N = 889$$

for all $k = 1, 2, \ldots$, and all $i = 0, 1, 2, \ldots$.

A corollary of this theorem is that there are infinitely many integers $m$ such $F_n + F_{n+1} + \ldots + F_{n+m} \equiv m \bmod N$, where $F_{n+i}$, $i = 0, 1, \ldots, m$, is the $(n+i)$-th Fermat number, and $n, N$ satisfying the hypothesis of the theorem.

In the second theorem of this paper (Theorem 2), the following result is proven. If $N = 2^p - 1$ and $p = 2^q - 1$ are both Mersenne primes, and $2^1 + 2^2 + 2^{2^2} + \dots + 2^{2^{q-1}}$ is not divisible by $N$, then the integers

$$1 + 2^{2^{kq}} + 2^{2^{kq+1}} + \dots + 2^{2^{kq+[(q-1)+rq+N_i]}}$$

are all divisible by $N$, for $k = 1, 2, \dots$ and $i = 0, 1, \dots$, where $r$ positive is some integer effectively found. (It satisfies a certain congruence; refer to Theorem 2.)

As a corollary to this, one sees that there are infinitely many integers m with the property that

$$F_{kq} + F_{kq+1} + \dots + F_{kq+m} \equiv (q-1) + rq \bmod (2^p - 1) \text{ for } k = 1, 2, \dots$$

**THEOREM 1.** Let $n$, $N$ be positive integers, $n, N \geq 2$, such that $2^n \equiv 1 \bmod s$, where s is the order of exponent to which 2 belongs modulo $N$. (Note that both $N$ and $s$ must be odd.) Let $l$, denote the largest integer such that $2^l < s$ (so $2^{l+1} \geq s$). Also let $b$ stand for the nonnegative integer so that $b < s$ and $2^{l+1} \equiv b \bmod s$. Set $k_l = 1 + 2 + 2^{2^1} + 2^{2^2} + \dots + 2^{2^l}$. If $a$ denotes the order or exponent to which 2 belongs mod $s$ ($s$ is odd and so $(2, s) = 1$), set

$$k_b = 2^b + 2^{b.2} + 2^{b.2^2} + \dots + 2^{b.2^{a-1}}.$$

Assume that $d \mid k_b$, where $d = (k_b, N)$.

Then the integers of the form $I_{n,m} = 1 + 2^{2^n} + 2^{2^{n+1}} + \dots + 2^{2^{n+m}}$ are all divisible by $N$ for infinitely many positive integers $m$. In fact, if $r$ is a positive integer such that $k_l + r.k_b \equiv 0 \bmod N$ (such an $r$ is guaranteed to exist on account of $d \mid k_l$ and it can be (effectively) computed), and $m = l + r.a + i.L.D$, where $L$ is the least common multiple of $N$ and $a$, $D = (a, N)$ and $i = 0, 1, 2, \dots$. Then $I_{n,m} \equiv 0 \bmod N$.

We present two examples.

1. Take $N = 31$, $n$ any positive integer with $n \equiv 0 \bmod 4$, say $n = 4k$. Then of course

$$2^n = (2^4)^k \equiv 1 \bmod 5, \text{ and } s = 5 \text{ (the exponent to which 2 belongs mod 31 is}$$

5). We have, in the notation of the theorem $l = 2$, and $2^{l+1} = 2^3 \equiv 3 \mod 5$, so $b = 3$. It is $a = 4$ and so $k_b = 2^3 + 2^{3.2} + 2^{3.2^2} + 2^{3.2^3}$, $k_l = 1 + 2^1 + 2^2 + 2^{2^2} = 23$.

We have, $3 \equiv 3, 3.2 \equiv 1, 3.2^2 \equiv 2$ and $3.2^3 \equiv 4 \mod 5$. Thus,

$k_b = 2^3 + 2 + 2^2 + 2^4 = 30 \mod 31$, which shows that $d = (k_b, N) = 1$.

Thus the hypothesis $d \mid k_l$ is satisfied.

To find an $r$ such that $k_l + r.k_b \equiv 0 \mod 31$, we must solve the congruence $23 + 30r \equiv 0 \mod 31$. The solution is $r \equiv 23 \mod 31$, and so the smallest positive solution is $r = 23$. Thus if we take $m = l + r.a + i.L.D. = 2 + 23.4 + 4.31i = 94 + 124i$ and since $n = 4k$, we conclude by the theorem that

$$I_{n,m} = 1 + 2^{2^{4k}} + 2^{2^{4k+1}} + \ldots + 2^{2^{4k+94+124i}} \equiv 0 \mod 31,$$

for all $k = 1, 2, \ldots$, and $i = 0, 1, 2, \ldots$.

2. Take $N = 7.31 = 217$, $n$ any positive integer with $n \equiv 0 \mod 4$, $n = 4k$. It is $2^{4k} = (2^4)^k \equiv 1 \mod 15$, which shows that $2^n \equiv 1 \mod 15$. The exponent to which 2 belongs modulo 217, is the least common multiple of the exponents modulo 7 and 31, so $s = 3.5 = 15$. We have $l = 3$ and $2^{l+1} = 2^4 \equiv 1 \mod 15$ which shows that $b = 1$. Also $a$ = least common multiple of the orders of 2 modulo 3 and 5 = 4. It is

$$k_l = 1 + 2^1 + 2^2 + 2^{2^2} + 2^{2^3} \text{ and } k_b = 2^1 + 2^2 + 2^{2^2} + 2^{2^3}$$

$$= 1 + 2 + 4 + 16 + 256 \qquad\qquad = 2 + 4 + 16 + 256$$

$$= 279 \qquad\qquad\qquad\qquad\qquad = 278.$$

Then $d = (278, 217) = 1$ and so $d \mid k_l$. We must solve the congruence

$279 + r.278 \equiv \mod (7.31)$. This congruence, when viewed mod 7, gives $6 + 5r \equiv 0 \mod 7$ which has a solution $r \equiv 3 \mod 7$. When considered mod 31, the above congruence gives, since $279 \equiv 0 \mod 31$, $278r \equiv 0 \mod 31$ which has solution $r \equiv 0 \mod 31$. The simultaneous congruences $r \equiv 3 \mod 7$ and $r \equiv 0 \mod 31$ have a solution $r \equiv 31 \mod 217$ (Chinese remainder theorem). The smallest positive integer with this property is 31. If we therefore take $m = l + r.a + i.L.D = 3 + 31.4 + 4.217i = 3 + 124 + 868i = 127 + 868i$, $n = 15k$, we conclude by the theorem that

$$I_{n,m} = 1 + 2^{2^{4k}} + 2^{2^{4k+1}} + \ldots\ldots + 2^{2^{4k+127+868i}} \equiv 0 \mod 217$$

for all $k = 1, 2, \ldots$, and $i = 0, 1, 2, \ldots$.

*Proof.* Since $2^n \equiv 1 \mod s$, we may set $2^n = k.s + 1$, for some positive integer $k$. Consider the integer $2^{2^n}$. We have $2^{2^n} = 2^{ks+1}$ which implies

$$2^{2n} \equiv 2^{ks+1} \equiv (2^s)^k .2 \equiv 1.2 \equiv 2 \mod N. \tag{1}$$

It immediately follows from congruence (1) that

$$2^{2^{n+j}} \equiv (2^{2^n})^{2^j} \equiv 2^{2^j} \mod N, \text{ for } j = 0, 1, 2, \ldots\ldots \tag{2}$$

If $k_l = 1 + 2^1 + 2^{2^1} + 2^{2^2} + \ldots + 2^{2^l}$, congruence (2) implies

$$2^{2^n} + 2^{2^{n+1}} + \ldots + 2^{2^{n+l}} = 2 + 2^{2^1} + \ldots + 2^{2^l} \equiv k_l \mod N. \tag{3}$$

According to the hypothesis, $l$ is the largest integer such that $2^l < s$, and so $2^{l+1} \equiv b \mod s$, for some $0 \leq b < s$. Since $a$ is the exponent to which 2 belongs mod $s$, it follows that $2^{a i} \equiv 1 \mod s$, $i = 0,1,2,...$ Given an $i$, we may set

$$2^{a i} = k_i.s + 1, \text{ for some } k_i \in \mathbb{Z}^+. \quad (4)$$

Also, since $s$ is the exponent to which 2 belongs modulo $N$, we have $2^s \equiv 1 \mod N$. Thus

$$2^{k_i s} \equiv \left(2^s\right)^{k_i} \equiv 1 \mod N \quad (5)$$

Congruences (4) and (5) imply

$$2^{b.2^{j}+i a} \equiv \left(2^{b.2^j}\right)^{2^{i a}} \equiv 2^{b.2^j}.\left(2^{k_i s}\right)^{b.2^j} \equiv 2^{b.2^j} \mod N \quad (6)$$

for $0 \leq j \leq a - 1$ and $i = 0,1,2,...$

If $k_b = 2^b + 2^{2.b} + 2^{b.2^2} + ... + 2^{b.2^{a-1}}$, it follows by (2), (3), (6) and

$$k_l = 1 + 2^1 + 2^{2^1} + 2^{2^2} + ... + 2^{2^l} \text{ that}$$

$$2^{2^n} + 2^{2^{n+1}} + ........ + 2^{2^{n+l+i a}} \equiv k_l + i k_b \mod N, \quad \text{for} \quad i = 0,1,2,... \quad (7)$$

By the hypothesis of the theorem, we have $d \mid k_b$, where $d = (k_b, N)$. Consequently there is a positive integer $r$ such that

$$k_l + r k_b \equiv 0 \mod N \quad (8)$$

Furthermore, if $L$ is the least common multiple of $N$ and $a$ and $D = (N, a)$, we must have

$$\frac{D.L}{a} \equiv 0 \mod N \quad (9)$$

Hence, if we replace $i$ by $r + i \frac{D.L}{a}$ in (7), we obtain

$$2^{2^n} + 2^{2^{n+1}} + \ldots + 2^{2^{n+1+r \cdot a + i.D.L}} \equiv k_l + \left( r + \frac{i.D.L}{a} \right) k_b$$

$$\equiv k_l + r k_b + \frac{i.D.L}{a} k_b \equiv 0 \bmod N,$$

by virtue of (8) and (9), and the theorem is proven. ∎

COROLLARY 1. *Let $n$, $N$ be integers satisfying the hypothesis of Theorem 1, with $N$ a composite integer and $2^n \geq N$. Then if $a$ is any proper positive divisor of $N$, there are infinitely many integers m such that the integers $\frac{1}{d}\left(1 + 2^{2^n} + 2^{2^{n+1}} + \ldots + 2^{2^{n+m}}\right)$ are composite integers divisible by $\frac{N}{d}$.*

*Proof.* By Theorem 1, there are infinitely many integers $m$ such that

$$I_{n,m} = 1 + 2^{2^n} + 2^{2^{n+1}} + \ldots + 2^{2^{n+m}} \equiv 0 \bmod N.$$

Hence $\frac{1}{d}\left(1 + 2^{2^n} + 2^{2^{n+1}} + \ldots + 2^{2^{n+m}}\right) \equiv 0 \bmod \frac{N}{d}$ for infinitely many $m$. Of course $1 < \frac{N}{d} < N$, since $d$ is a proper divisor of $N$. From $2^n \geq N$ it follows that $1 + 2^{2^n} \geq 1 + 2^N > N$, which implies

$$1 + 2^{2^n} + 2^{2^{n+1}} + \ldots + 2^{2^{n+m}} > N,$$

whence

$$\frac{1}{d}\left(1+2^{2^n}+2^{2^{n+1}}+\ldots+2^{2^{n+m}}\right)=\frac{I_{n,m}}{d}>\frac{N}{d}$$

which together with $\frac{I_{n,m}}{d}\equiv 0\bmod\left(\frac{N}{d}\right)$, establishes that the integers $\frac{I_{n,m}}{d}$ are composite for infinitely many integer values of *m*. ∎

COROLLARY 2. Let *n, N* be integers satisfying the hypothesis of Theorem 1. Then there are infinitely many integers *m* so that $F_n+F_{n+1}+\ldots+F_{n+m}\equiv m\bmod N$, where $F_{n+i}$, $i=0,1,2,\ldots,m$, stands for the $(n+i)$-th Fermat number.

*Proof.* According to Theorem 1, there are infinitely many numbers *m* so that

$$1+2^{2^n}+2^{2^{n+1}}+\ldots+2^{2^{n+m}}\equiv 0\bmod N$$

We have $F_{n+i}=2^{2^{n+i}}+1$, for $i=0,1,2,\ldots,m$. Thus,

$$F_n+F_{n+1}+\ldots+F_{n+m}=\left(2^{2^n}+1\right)+\left(2^{2^{n+1}}+1\right)+\ldots\ldots+\left(2^{2^{n+m}}+1\right)$$

$$=\left(1+2^{2^n}+2^{2^{n+1}}+\ldots+2^{2^{n+m}}\right)+[(m+1)-1]$$

= m mod N,

by virtue of (10). ∎

THEOREM 2. Let $N = 2^P - 1$ and $p = 2^q - 1$ be both Mersenne

primes (such

is the case for instance when *q* = 3 and *p* = 7; *N* = 127). Then

$$1+2^{2^{kq}}+2^{2^{kq+1}}+\ldots+2^{2^{kq+[(q-1)+rq+Ni]+m}}\equiv 0\bmod N$$

for all *k* = 1,2,..., *i* = 0,1,2,..., where *r* is the smallest positive integer such that

$$\left(1+2^1 + 2^{2^1} + 2^{2^2} \ldots + 2^{2^{q-1}}\right)+(2^1 + 2^{2^1} + 2^{2^2} \ldots + 2^{2^{q-1}})r \equiv 0 \bmod N$$

provided that $2^1 + 2^{2^1} + 2^{2^2} \ldots + 2^{2^{q-1}}$ is no divisible by $N$

*Proof.* Let $n=q.k$, $k$ being any positive integer. It is $2^n = \left(2^q\right)^k = (p+1)^k \equiv 1 \bmod p$. However, it is clear that $s = p$, where $s$ is the exponent to which 2 belongs modulo $N$. This follows from the fact that $2^p = N + 1$ and $p$ is prime. (Thus the exponent $x$ to which 2 belongs modulo $N$ would have to be a divisor of $p$, since $2^p \equiv 1 \bmod N$. Therefore, $x = 1$ or $p$, since $p$ is a prime. It cannot be $x = 1$, since $N > 1$, thus $x = p$.)

We have thus proved $2^n \equiv 1 \bmod s$, which means that Theorem 1, since $n, N \geq 2$, does apply. In the notation of that theorem, we have $l = q - 1$ (This is because $2^{q-1} < s = p = 2^q - 1$, in virtue of $q \geq 3$, in particular $q \geq 2$; and because $2^{l+1} = 2^q > s = p = 2^q - 1$) The exponent $a$ to which 2 belongs modulo $s = p$, equals $q$, $a = q$. (The reasoning is the same by showing that $p$ is the exponent to which 2 belongs modulo $N$, which was done above.)

Also, $2^{l+1} = 2^q = p+1 \equiv 1 \bmod p$ which shows that $b = 1$. Thus

$$k_l = 1 + 2^1 + 2^{2^1} + 2^{2^2} + \ldots + 2^{2^l} = 1 + 2^1 + 2^{2^1} + 2^{2^2} + \ldots + 2^{2^{q-1}}$$

and

$$k_b = 2^1 + 2^{1.2} + 2^{1.2^2} + \ldots + 2^{2^{q-1}}.$$

If, as it is assumed in the hypothesis, $N$ does not divide $k_b$, then, since $N$ is a prime, it follows that $d = (k_b, N) = 1$. Thus $d | k_l$, which shows that the hypothesis of Theorem 1 is satisfied. (Note that if $k_b \equiv 0 \bmod N$, then, since $(k_l, k_b) = 1$ and by virtue of $k_l = 1 + k_b$, the hypothesis of Theorem 1 cannot be satisfied). Consequently there is a smallest positive integer $r$ satisfying $k_l + rk_b \equiv 0 \mod N$. By taking $m = l + r.a + i.L.D$, where $L$ = the least common multiple of $N$ and $a = (2^p - 1)q$ (since $2^p - 1$ is a prime and $q < 2^p - 1$, by virtue of $p = 2^q - 1$, it follows that $(2^p-1, q) = 1$; thus $L = q(2^p-1)$) and $D = (N, a) = (2^p - 1, q) = 1$. We find that $m = q - 1 + r.q + i.q.(2^p - 1)$, and therefore, by applying Theorem 1, the proof is concluded. ∎

COROLLARY 3. If $N, p$ are Mersenne primes, satisfying the hypothesis of Theorem 2, i.e., $N = 2^p - 1$, $p = 2^q - 1$, $F_{qk} + F_{qk+1} + \ldots + F_{qk[q-1+r.q+i.q.(2^p-1)]} \equiv q - 1 + rq \mod(2^p - 1)$ for all $k = 1, 2, \ldots,$ and $i = 0, 1, 2, \ldots$.

*Proof.* It follows immediately by Corollary 2 to Theorem 1, and by Theorem 2, and by virtue of $m = q - 1 + rq + iq.(2^p - 1)$. ∎